%% file: bcv.tex
\documentclass{amsart}
\usepackage{amscd,graphicx}

\include{loopmacro}

\input xy
\xyoption{all}

\begin{document}

\setcounter{page}{1}

\title[Yet another delooping machine]{Yet another delooping machine}

\author[B. Badzioch]{Bernard Badzioch}
\author[K. Chung]{Kuerak Chung}
\author[A. A. Voronov]{Alexander A. Voronov}

\address {School of Mathematics\\University of Minnesota\\
Minneapolis, MN 55455}
\email{badzioch@math.umn.edu}
\email{krchung@math.umn.edu}
\email{voronov@math.umn.edu}

\date{March 4, 2004}

\thanks{The third author was supported in part by NSF grant DMS-0227974}

\begin{abstract}
  We suggest a new delooping machine, which is based on recognizing an
  $n$-fold loop space by a collection of operations acting on it, like
  the traditional delooping machines of Stasheff, May, Boardman-Vogt,
  Segal, and Bousfield. Unlike in the traditional delooping machines,
  which carefully select a nice space of such operations, we consider
  all natural operations on $n$-fold loop spaces, resulting in the
  algebraic theory $\Map (\bigvee_\bullet S^n, \bigvee_\bullet S^n)$.
  The advantage of this new approach is that the delooping machine is
  universal in a certain sense, the proof of the recognition principle
  is more conceptual, works the same way for all values of $n$, and
  does not need the test space to be connected.
\end{abstract}

\maketitle

\section{Introduction}

\label{INTRO}

\begin{figure}
\label{figure}
\centerline{\includegraphics[width=3.2in]{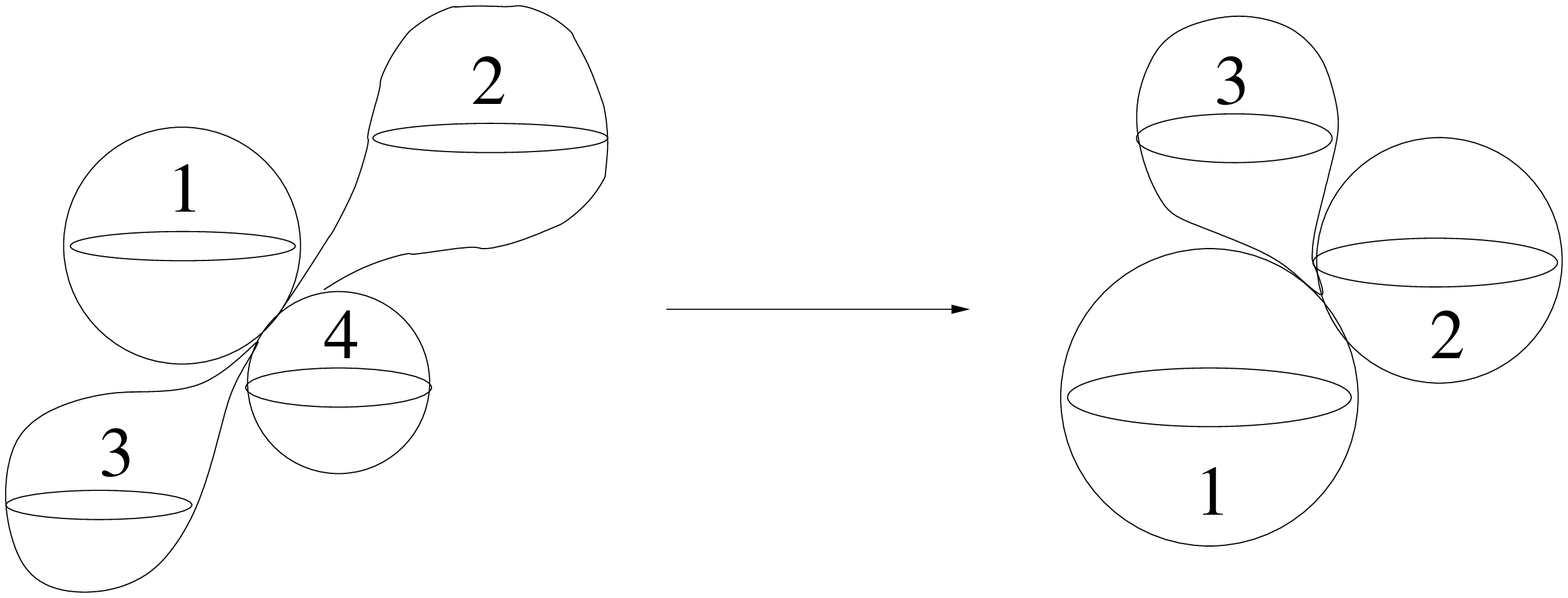}}
\caption{A morphism $\bigvee_3 S^n \to \bigvee_4 S^n$ in category $\Tn$.}
\label{k=1}
\end{figure}

The goal of this paper is to give a proof of the following
characterization of $n$--fold loop spaces. In the category $\Spaces$
of pointed spaces, consider the full subcategory generated by the
wedges $\bigvee_k S^n$ of $n$-dimensional spheres for $k\geq 0$ (where
$\bigvee_0 S^n = \ast$). Let $\Tn$ denote the opposite category, see
Figure~\ref{figure}. Since $\bigvee_k S^n$ is a $k$-fold coproduct of
$S^n$'s in $\Spaces$, in $\Tn$ it is a $k$-fold categorical product of
$S^n$'s.

\begin{thm}
\label{MAIN}
A space $Y\in \Spaces$ is weakly equivalent to an $n$--fold loop
space, iff there exists a product preserving functor $\Xtilde \colon
\Tn \ra \Spaces$ such that $\Xtilde (S^n)\simeq Y$.
\end{thm}

The category $\Tn$ is in fact an \emph{algebraic theory} (see
\ref{STRICT ALG}). From this point of view, one can regard the above
theorem as a \emph{recognition principle}: a loop space structure is
detected by the structure of an \emph{algebra over the algebraic
  theory} $\Tn$.

We will actually prove a stronger version (see Theorem~\ref{MAIN2}) of
Theorem~\ref{MAIN}: given a product preserving functor $\Xtilde: \Tn
\to \Spaces$, one can construct a space $B_n \Xtilde$ such that
$\Omega^n B_n \Xtilde \simeq \Xtilde (S^n)$, thereby \emph{delooping}
the space $\Xtilde (S^n)$.

This description of iterated loop spaces is in some sense an extreme
delooping machine. By Yoneda's lemma the theory $\Tn$ encodes all
natural maps $(\Omega^n X)^k \ra (\Omega^n X)^l$, and we use all this
structure in order to detect loop spaces.  This stands in contrast to
the approach of Stasheff \cite{stasheff}, May \cite{May},
Boardman-Vogt \cite{boardman-vogt}, Segal \cite{segal}, or Bousfield
\cite{Bousfield}, where only carefully chosen sets of maps of of loop
spaces are used for the same purpose.  Our indiscriminate method
however brings some advantages.  First of all, as in \cite{Bousfield},
Theorem \ref{MAIN} is true for all, not necessarily connected, loop
spaces.  Also, since we avoid making particular choices of operations
on loop spaces, thus constructed delooping machine provides a
convenient ground for proving uniqueness theorems of the kind of May
and Thomason \cite{may-thomason}, \cite{thomason}. Namely, given an
operad, a PROP, or a semi-theory (i.e., a machine of the type of
Segal's $\Gamma$-spaces, see \cite{badzioch:gamma}), one can replace 
it by an algebraic theory describing the same structure on spaces. On the other
hand, it is relatively easy to compare homotopy theories of objects
described by various algebraic theories. This implies Theorem~\ref{unique} --
a uniqueness result  for "delooping theories".

Most of the arguments and constructions we use are formal and do not
depend on any special properties of loop spaces. Indeed, at least one
implication of the statement of Theorem \ref{MAIN} holds when we
replace $S^n$ with an arbitrary pointed space $A$. If $\TT^A$ is an
algebraic theory constructed analogously to $\Tn$ above, then for any
mapping space $Z= \Map(A, Y)$, we can define a product preserving
functor $\Xtilde\colon \TT^A \ra \Spaces$ such that $\Xtilde(A) = Z$.
We do not expect that for an arbitrary $A$ also the opposite statement
will be true, that is that any such functor will come from some space
$ \Map(A, Y)$.  It should be true, however, that if for a given space
$A$, the mapping spaces from $A$ can be described as algebras over
some operad, PROP, semi-theory, algebraic theory, or using some other
formalism employing only finitary operations on a space, then they
must be characterized by means of the theory $\TT^A$.

Another advantage of the proposed recognition principle is that the
argument seems to be more conceptual than in the previously known
cases. For example, we get an analogue (Corollary~\ref{APX}) of May's
approximation theorem \cite{May} as a simple consequence of, rather
than a hard step towards the recognition principle.

This simplicity comes, no wonder, with a price tag attached: while the
homology of the little $n$-disks operad has a neat description as the
operad describing $n$-algebras, see F.~Cohen \cite{cohen:cubes,
  cohen:sundry}, even the rational homology of the corresponding PROP
$\Map (\bigvee_l S^n, \bigvee_k S^n)$ is harder to come by, see the
thesis \cite{chung} of the second author.

The theory $\Tn$ bears resemblance to the cacti operad of
\cite{sullivan-v}, which consists of (unpointed) continuous maps from
a sphere $S^n$ to a tree-like joint of spheres $S^n$ at finitely many
points. This operad was invented as a bookkeeping device for
operations on free sphere spaces arising in string topology, see
\cite{chas-sullivan}.

Also, the operadic part $\OO_n := \Map (S^n, \bigvee S^n)$ of  
$\Tn$ has been described as a "universal operad of $n$-fold
loop spaces'' by P.~Salvatore in \cite{salvatore}. As it was also
noted by Salvatore, while the space underlying an algebra over this
operad is weakly equivalent to an $n$-fold loop space, in general a
loop space will admit several actions of $\OO_n$. Therefore
$\OO_n$-algebras can be seen as loop spaces equipped with some extra
structure.

\begin{notation}
  \hfill

\begin{itemize}
\item Let $\Spaces$ denote the category of pointed compactly generated
  (but not necessarily Hausdorff) topological spaces. From the
  perspective of homotopy theory, there is no difference between this
  category and the category of all pointed topological spaces.  The
  category $\Spaces$ has a model category structure with the usual
  notions of weak equivalences, fibrations and cofibrations, and it is
  Quillen equivalent to the category of pointed topological spaces,
  see \cite{hovey}. The assumption that all spaces are compactly
  generated has the advantage that for any space $X$, the smash
  product functor $Y \mapsto Y\wedge X$ is left adjoint to the mapping
  space functor $Y \mapsto \Map(X, Y)$.  This has some further useful
  consequences which we will invoke.
  
\item If $X$ is an unpointed space by $X_+$ we will denote the space
  $X$ with an adjoined basepoint.

\item  All functors are assumed to be covariant.
  
\item If $\CC$ is a category, then $\CC^\op$ will denote the opposite
  category of $\CC$.

\end{itemize}
\end{notation}

\section{Algebraic theories and their algebras} 

\label{THEORIES}

\begin{df}
\label{STRICT ALG}

An \emph{algebraic theory} $ \TT$ is a category with objects $T_0,
T_1, \dots$ together with, for each $n$,  a choice of morphims
$p^n_1, \dots, p^n_n \in \Mor_\TT(T_n, T_1)$ such that for any $k, n$
the map 
$$\prod_{i=1}^n p^n_i \colon \Mor_\TT(T_k, T_n)\ra
\prod_{i=1}^{n}\Mor_\TT(T_k, T_1)$$
is an isomorphism. In other words,
the object $T_n$ is an $n$-fold categorical product of $T_1$'s, and
$p^n_i$'s are the projection maps. In particular $T_0$ is the terminal
object in $\TT$. We will also assume that it is an initial object. A
morphism of algebraic theories is a functor $\TT\ra \TT'$ preserving
the projection maps. We will consider algebraic theories enriched over
$\Spaces$; in particular, the sets of morphisms will be provided with
a pointed topological space structure.

Given an algebraic theory $\TT$, a $\TT$-\emph{algebra} $\Xtilde$ is a
product preserving functor $\Xtilde\colon\TT\ra \Spaces$. A morphism
of $\TT$-algebras is a natural transformation of functors.
\end{df}

We will say that a space $X$ admits a $\TT$-algebra structure, if
there is a $\TT$-algebra $\Xtilde$ and a homeomorphism $\Xtilde(T_1)\cong X$.

For an algebraic theory $\TT$, by $\Alg_\TT$ we will denote the category of 
all $\TT$-algebras and their morphisms.

\begin{ex}
\label{T_A}
For any pointed space $A\in \Spaces$ we can define an algebraic theory
$\TT_A$ enriched over $\Spaces$ by setting
$$\Mor_{\TT_A}(T_m, T_n) := \Map(A^m, A^n).$$
Thus, $\TT_A$ is
isomorphic to the full subcategory of $\Spaces$ generated by the
spaces $A^n$ for $n\geq 0$.

For any $Y\in \Spaces$, we can consider a product preserving functor
$$\TT_A\ra \Spaces, \hskip 1cm T_n\mapsto \Map(Y, A^n).$$
This shows
that any mapping space $\Map(Y, A)$ has a canonical structure of a
$\TT_A$-algebra.
\end{ex}

\begin{ex}
\label{T^A}
Let $A$ be again a pointed space,  and let $\TT^A$ be a category with objects 
$T_0, T_1, \dots$ and morphisms
$$\Mor_{\TT^A}(T_m, T_n)=\Map(\bigvee_n A, \bigvee_m A).$$
In other
words, $\TT^A$ is isomorphic to the opposite of the full subcategory
of $\Spaces$ generated by the finite wedges of $A$. Since $\bigvee_n
A$ is an $n$-fold coproduct of $A$ in $\Spaces$, $T_n$ is an $n$-fold
categorical product of $T_1$'s in $\TT^A$. It follows that $\TT^A$ is
an algebraic theory. For $Y\in \Spaces$, we can define a functor
$$\TT^A\ra \Spaces, \hskip 1cm T_n \mapsto \Map(\bigvee_n A, Y).$$
Therefore the mapping space $\Map(A, Y)$ has a canonical structure of
a $\TT^A$-algebra.  In particular, if $A=S^n$ we get that any $n$-fold
loop space canonically defines an algebra over $\Tn$.
\end{ex}

\begin{nn}
\label{S^0}
A special instance of an algebraic theory $\TT^A$ is obtained when we
take $A= S^0$.  The category ${\TT^{S^0}}$ is equivalent to the
opposite of the category of finite pointed sets. One can check
that the forgetful functor
$$U_{\TT^{S^0}} \colon \Alg^{\TT^{S^0}} \ra \Spaces, \hskip 1cm
U_{\TT^{S^0}}(\Xtilde)=\Xtilde(T_1),$$
gives an equivalence of
categories. Also, for any algebraic theory $\TT$ there is a unique map
of algebraic theories $I_\TT\colon \TT^{S^0}\ra \TT$.  
If $U_\TT\colon \Alg^\TT \ra \Spaces$ is the forgetful functor, $U_\TT(\Xtilde)
= \Xtilde(T_1)$, then we have $U_\TT = U_{\TT^{S^0}}\circ
{I_\TT}^\ast$ where ${I_T}^\ast\colon \Alg^\TT \ra \Alg^{\TT^{S^0}}$
is the functor induced by $I_\TT$.
\end{nn}

\section{Tensor product of functors} 


\begin{df}
  Let $\CC$ be a small topological category, i.e., a small category
  enriched over $\Spaces$, and $F\in \Spaces^\CC$, $G\in
  \Spaces^{\CC^{\rm op}}$. The tensor product $F\otimes_\CC G$ is the
  colimit
$$\xymatrix{
F\otimes_\CC G  := \colim  \bigvee_{(c, d)\in \CC\times \CC} \Mor(c, d)\wedge 
F(c)\wedge G(d) \ar@<0.5ex>[r]^-{j_1} 
\ar@<-0.5ex>[r]_-{j_2} &  
\bigvee_{c \in \CC} F(c)\wedge G(c)
}.$$
The map $j_1$ is the wedge of the maps 
$\ev\wedge \id \colon (\Mor(c, d)\wedge F(c))\wedge G(d) \ra
F(d)\wedge G(d)$,
where $\ev$ is the evaluation map, and $j_2$ is similarly induced by
the evaluation maps
$\ev \colon \Mor(c, d)\wedge G(d) \ra G(c)$.
\end{df}  

The most important -- from our perspective -- property of the tensor product 
is given by the following
\begin{prop} 
\label{ADJ TENSOR}
Let $\CC$ be a small topological category and $G\in
\Spaces^{\CC^{\rm op}}$.  Consider the functor
$$\Map(G, -)\colon  \Spaces\ra \Spaces^\CC, \hskip 1cm  Z\mapsto \Map(G, Z).$$
The left adjoint of $\Map(G, - )$ exists and is given by 
$$-\otimes_\CC G \colon \Spaces^\CC \ra \Spaces, \hskip 1cm  F\mapsto F\otimes_\CC G.$$
\end{prop}
\noindent For a proof see, e.g., \cite{ML}.




\begin{nn}
\label{TENSOR FUNCTORS}

Assume now that we have two small categories $\CC$ and $\DD$ enriched
over $\Spaces$ and two functors $F\colon \CC\times \DD \ra \Spaces$
and $G\colon \CC^\op \ra \Spaces$. For every $d\in \DD$, the functor
$F$ defines $F(d)\colon \CC \ra \Spaces$ by $F(d)(c) = F(c, d)$.
Applying the tensor product construction, we obtain a new functor
$F\otimes_\CC G \colon \DD \ra \Spaces$ such that $(F\otimes_\CC G)(d)
= F(d)\otimes_\CC G$. Since smash product in $\Spaces$ commutes with
taking colimits, for any $H\colon \DD^\op\ra \Spaces$ we have a
natural isomorphism
$$H\otimes_{\DD^\op}(F\otimes_\CC G)\cong
(H\otimes_{\DD^\op}F)\otimes_\CC G \in \Spaces.$$
\end{nn}

\begin{nn}
\label{EXAMPLES OF TENSORS}

Our main interest lies in the following instances of these constructions:

1) For $A\in \Spaces$, let $\TT^A$ be the algebraic theory defined in
Example~\ref{T^A}.  Consider the functor
$$\Omega^A\colon \Spaces \ra \Spaces^{\TT^A} $$
given by
$\Omega^A(Y)(T_k) := \Map(\bigvee_k A, Y)$. By Proposition \ref{ADJ
  TENSOR}, $\Omega^A$ has a left adjoint $B_A\colon \Spaces^{\TT^A}\ra
\Spaces$, given by $B_A(F)= F\otimes_{\TT^A} \bigvee_\bullet A$. Here
$\bigvee_\bullet A$ denotes the functor from $(\TT^A)^\op$ to
$\Spaces$ such that $\bigvee_\bullet A(T_k) = \bigvee_k A$. Note that
$\Omega^A (Y)$ preserves products, and so $\Omega^A$ takes values in
the full subcategory $\Alg^{\TT^A} \subset \Spaces^{\TT^A}$.  Thus, we
get an adjoint pair $(B_A, \Omega^A)$ of functors between
$\Alg^{\TT^A}$ and $\Spaces$.

2)  For  $A\in \Spaces$,  let $\End{\bigvee_\bullet A}$ denote the 
functor $\TT^A \times (\TT^{A})^\op \to \Spaces$ defined by 
$$\End{\bigvee_\bullet A}(T_k, T_l):= \Map(\bigvee_k A, \bigvee_l
A).$$
Using the canonical map $I_{\TT^A}\colon \TT^{S^0}\ra \TT^A$
(see \ref{S^0}), we can view $\End{\bigvee_\bullet A}$ as a functor on
the category $\TT^A \times (\TT^{S^0})^\op$.  For $Y\in \Spaces$,
define
$$F_{\TT^A}(Y):=\End{\bigvee_\bullet A}\otimes_{(\TT^{S^0})^\op}
\Omega^{S^0}(Y) \in \Spaces^{\TT^A}.$$
One can check that
$F_{\TT^A}(Y)$ preserves products, i.e., defines a $\TT^A$-algebra.
Thus we get a functor
$$F_{\TT^A}\colon \Spaces \ra \Alg^{\TT^A}, \   \   Y\mapsto  F_{\TT^A}(Y),$$
which is left adjoint to the forgetful functor 
$$U_{\TT^A}\colon \Alg^{\TT^A} \ra \Spaces, \ \ U_{\TT^A}(\Xtilde)=
\Xtilde(T_1).$$
We will call $F_{\TT^A}$ the free $\TT^A$-algebra
functor and $F_{\TT^A}(Y)$ the free $\TT^A$-algebra generated by $Y$.

3) Consider again an algebraic theory $\TT^A$ and let $\Delta^\op$ be
the simplicial category.  Let $\Xtilde_\bullet \colon \TT^A\times
\Delta^\op \ra \Spaces$ be a simplicial $\TT^A$--algebra. Let
$\Delta[\bullet]_+ \colon \Delta \ra \Spaces$ denote the pointed
cosimplicial space $[n]\mapsto \Delta[n]_+$.  In this case the tensor
product $\Xtilde_\bullet \otimes_{\Delta^\op} \Delta[\bullet]_+ =:
|\Xtilde_\bullet | $ gives the geometric realization of
$\Xtilde_\bullet$. Since realization preserves products in $\Spaces$,
we see that $|\Xtilde_\bullet |$ is a $\TT^A$--algebra.
\end{nn}

\begin{nn}
\label{FREE ALG}

Notice that the isomorphism of Section~\ref{TENSOR FUNCTORS} shows
that for a pointed simplicial space $Y_\bullet$ we have $|F_{\TT^A}
Y_\bullet|\cong F_{\TT^A}|Y_\bullet |$, and that similarly for a
simplicial $\TT^A$--algebra $\Xtilde_\bullet$ we get $|B_A
\Xtilde_\bullet | \cong B_A|\Xtilde_\bullet |$.
\end{nn}

\begin{nn}
\label{BF}

Finally, consider the functors $\Omega^A$ and $U_{\TT^A}$ of
Section~\ref{EXAMPLES OF TENSORS}. The composition $U_{\TT^A}\circ
\Omega^A\colon \Spaces \ra \Spaces$ is given by $U_{\TT^A}\circ
\Omega^A(Y) = \Map(A, Y)$. As a result its left adjoint $B_A \circ
F_{\TT^A}$ is the smash product $B_A \circ F_{\TT^A}(Y)= Y\wedge A$.
This observation indicates that the algebraic theory $\TT^A$ may be
suitable for describing mapping spaces from $A$, at least in some
cases. Indeed, a simple computation shows that for a finite pointed
set $Z$, we have $\Mor_{\Alg^{\TT^A}} (\Map(A, Z \wedge A), \Xtilde)
\cong \Map (Z, U_{\TT^A} (\Xtilde))$.  Thus, by the adjointness of
$F_{\TT^A}$ and $U_{\TT^A}$, we get

\begin{lm}
\label{KEY}
For any pointed finite set $Z$, we have a canonical isomorphism
\[
F_{\TT^A}Z \cong \Map(A, Z\wedge A)
\]
of $\TT^A$-algebras.
\end{lm}

\noindent Combining this isomorphism with the equality 
$B_A (F_{\TT^A}(Z))= Z\wedge A$, we see that $B_A$ acts as a
classifying space for $\Map(A, Z\wedge A)$. Our goal will be to show
that when we take $A=S^n$, this construction works for any $\Tn$-
algebra.
\end{nn}

\section{Model categories and Quillen equivalences} 

\label{MODELS}

Our strategy of approaching Theorem \ref{MAIN} will be to reformulate
it in the language of model categories and prove it in this form.
Below we describe model category structures we will encounter in this
process. As it was the case so far, most of our setup will apply to
mapping spaces $\Map(A, Y)$ from an arbitrary space $A$, and only in
the proof of Theorem~\ref{MAIN2}, we will specialize to $A=S^n$.

For any algebraic theory $\TT$, the category of $\TT$-algebras
$\Alg^\TT$ has a model category structure with weak equivalences and
fibrations defined objectwise, i.e., via the forgetful functor
$U_\TT$, \cite{schwanzl-vogt}.  For a CW-complex $A\in \Spaces$, let
$R_A\Spaces$ denote the category of pointed spaces together with the
following choices of classes of morphisms:
\begin{itemize}
\item[-] a map $f\colon Y\ra Z$ is a weak equivalence in $R_A\Spaces$,
  if $f_\ast\colon \Map(A, Y)\ra \Map(A,Z)$ is a weak equivalence of
  mapping spaces;
\item[-] a map $f$ is a fibration if it is a Serre fibration;
\item[-] a map $f$ is a cofibration if it has the left lifting
  property with respect to all fibrations which are weak equivalences
  in $R_A\Spaces$.
\end{itemize}

\begin{prop}
\label{R_A}
The category $R_A\Spaces$ is a model category.
\end{prop}
  
\begin{proof}
  The statement follows from a general result on the existence of
  right localizations of model categories, see \cite[5.1, p.
  65]{hirschhorn}.
\end{proof}
Note that for $A = S^0$, this defines the standard model category
structure on $\Spaces$.

In order to avoid confusing $R_A\Spaces$ with $\Spaces$, we will call
weak equivalences (respectively, fibrations and cofibrations) in
$R_A\Spaces$ \emph{$A$-local equivalences} (respectively,
\emph{fibrations} and \emph{cofibrations}). Notice that a map $f\colon
Y\ra Z$ is an $S^n$-local equivalence, iff it induces isomorphisms
$f_\ast \colon \pi_q(Y)\ra \pi_q(Z)$ for $q\geq n$.

\begin{nn} \textbf{A cofibrant resolution of a $\TT^A$-algebra}.
\label{FUA}
Directly from the definition of the model structure on $\Alg^{\TT^A}$,
it follows that every $\TT^A$-algebra is a fibrant object. The
structure of cofibrant algebras is more complicated (see
\cite{schwanzl-vogt}).  For an arbitrary algebra $\Xtilde\in
\Alg^{\TT^A}$, one can however describe its cofibrant replacement as
follows. Recall the adjoint pair
$$\xymatrix{ F_{\TT^A} \colon \Spaces \ar@<0.5ex>[r]& \Alg^{\TT^A}
  \colon U_{\TT^A} \ar@<0.5ex>[l] }$$
of Section~\ref{EXAMPLES OF
  TENSORS}.2.
\end{nn}

\begin{prop}
\label{Q PAIR 1}
For any CW-complex $A\in \Spaces$, the functors 
$$\xymatrix{ F_{\TT^A}\colon \Spaces \ar@<0.5ex>[r]& \Alg^{\TT^A}
  \colon U_{\TT^A} \ar@<0.5ex>[l] }$$
form a Quillen pair.
\end{prop}

\begin{proof}
  The functor $U_{\TT^A}$ sends weak equivalences and fibrations in
  $\Alg^{\TT^A}$ to weak equivalences and fibrations in $\Spaces$,
  respectively, thus the conclusion follows.
\end{proof}

Next, consider the adjoint functors
$$\xymatrix{ | \cdot |\colon \SSets \ar@<0.5ex>[r]& \Spaces \colon
  \Sing \ar@<0.5ex>[l] }$$
between the categories of pointed spaces
and pointed simplicial sets, where $\Sing$ is the singularization
functor and $|\cdot |$ is geometric realization. We will denote by
$F'_{\TT^A}\colon \SSets\ra\Alg^{\TT^A}$ the composition of $|\cdot |$
and $F_{\TT^A}$, and by $U'_{\TT^A}\colon \Alg^{\TT^A}\ra \SSets$ the
functor obtained by composing $U_{\TT^A}$ with $\Sing$. The functors
$F'_{\TT^A}$, $U'_{\TT^A}$ form again a Quillen pair. Therefore for
any $\TT^A$-algebra $\Xtilde$, they define a simplicial object
${F'_{\TT^A}U'_{\TT^A}}_\bullet \Xtilde$ in the category
$\Alg^{\TT^A}$ which has the algebra $(F'_{\TT^A}U'_{\TT_A})^{(k+1)}
\Xtilde$ in its $k$-th simplicial dimension. Its face and degeneracy
maps are defined using the counit and the unit of adjunction,
respectively (compare \cite[Chapter 9]{May}).  Let
$|{F'_{\TT^A}U'_{\TT^A}}_\bullet \Xtilde|$ denote the objectwise
geometric realization of ${F'_{\TT^A}U'_{\TT^A}}_\bullet \Xtilde$.

\begin{lm}
\label{FUA ALG}
 $|{F'_{\TT^A}U'_{\TT^A}}_\bullet \Xtilde|$ is a $\TT^A$-algebra.
\end{lm}

\begin{proof}
  Clearly, $|{F'_{\TT^A}U'_{\TT^A}}_\bullet \Xtilde|$ is a functor from
  $\TT^A$ to $\Spaces$.  Also, since we are working in the category of
  compactly generated spaces, realization preserves products, and so
  $|{F'_{\TT^A}U'_{\TT^A}}_\bullet \Xtilde|$ is a $\TT^A$-algebra.
\end{proof}

Similarly to \cite[3.5, p.~903]{ATHT}, we get

\begin{lm}
\label{FUA WE}
For any $\Xtilde\in \Alg^{\TT^A}$ there is a canonical weak equivalence
$$|{F'_{\TT^A}U'_{\TT^A}}_\bullet \Xtilde|\ra \Xtilde .$$
\end{lm}

The above lemma remains to be true, if we replace the functors
$F'_{\TT^A}$ and $U'_{\TT^A}$ with $F_{\TT^A}$ and $U_{\TT^A}$,
respectively.  What we will use in the sequel (see Step 3 of the proof
of Theorem~\ref{MAIN2})  though is that the free algebras
$(F'_{\TT^A}U'_{\TT^A})_n \Xtilde$ are generated by spaces obtained as 
realizations of simplicial sets. 
The algebra $|{F'_{\TT^A}U'_{\TT^A}}_\bullet \Xtilde|$ can be taken as
a cofibrant replacement of $\Xtilde$, since we have

\begin{lm}
\label{FUA COF}
For any $\Xtilde\in \Alg^{\TT^A}$ the algebra
$|{F'_{\TT^A}U'_{\TT^A}}_\bullet \Xtilde|$ is a cofibrant object in
$\Alg^{\TT^A}$.
\end{lm}

\begin{proof}
  This is a consequence of \cite{schwanzl-vogt}, which describes the
  structure of cofibrant  objects in the model category $\Alg^\TT$.
\end{proof}

Next, let $A\in \Spaces$. Recall (Section~\ref{EXAMPLES OF TENSORS}.1)
that we have an adjoint  pair of functors $(B_A, \Omega^A)$. Moreover
the following holds:

\begin{prop}
\label{Q PAIR}
For any CW-complex $A\in \Spaces$, the functors
$$\xymatrix{ B_A\colon \Alg^{\TT^A}\ar@<0.5ex>[r]& R_A\Spaces \colon
  \Omega^A \ar@<0.5ex>[l] }$$
form a Quillen pair.
\end{prop}

\begin{proof}
  The functor $\Omega^A$ sends $A$-local equivalences and $A$-local
  fibrations to weak equivalences and fibrations in $\Alg^{\TT^A}$,
  respectively which yields the statement follows.
\end{proof}

Our main result, Theorem~\ref{MAIN}, can now be restated more
precisely as follows:
\begin{thm}
\label{MAIN2}
For $n\geq 0$ the Quillen pair
$$\xymatrix{ B_n\colon \Alg^{\Tn}\ar@<0.5ex>[r]& R_{S^n}\Spaces \colon
  \Omega^n \ar@<0.5ex>[l] },$$
where $B_n := B_{S^n}$ and $\Omega^n :=
\Omega^{S^n}$, is a Quillen equivalence. In particular, the two
functors induce an equivalence of the homotopy categories.
\end{thm}

\begin{crl}[Approximation theorem]
\label{APX}
For any CW-complex $X \in \Spaces$, the following $\Tn$-algebras are
weakly equivalent:
\[
F_n X \xrightarrow{\sim} \Omega^n \Sigma^n X ,
\]
where $F_n X$ denotes the free $\Tn$-algebra $F_{\Tn} X$ on $X$ and
$\Sigma^n X = S^n \wedge X$ is the reduced suspension. Moreover, these
equivalences establish an equivalence of monads $F_n \sim \Omega^n
\Sigma^n$ on the category of CW-complexes.
\end{crl}

Let us first deduce Theorem~\ref{MAIN} and Corollary~\ref{APX} from
Theorem~\ref{MAIN2}.

\begin{proof}[Proof of Theorem~$\ref{MAIN}$]
  Let $\Xtilde$ be any $\Tn$-algebra, and let
  $\Xtilde\stackrel{\sim}{\ra} \Xtilde_c$ be its cofibrant
  replacement. Like any other object in $R_{S^n} \Spaces$, $B_n
  \Xtilde_c$ is fibrant and therefore Theorem~\ref{MAIN2} implies that
  the adjoint $\Xtilde_c \to \Omega^n B_n \Xtilde_c$ of the identity
  isomorphism $B_n \Xtilde_c \xrightarrow{\simeq} B_n \Xtilde_c$ is a
  weak equivalence of $\Tn$-algebras. Therefore $\Xtilde(T_1)\simeq
  \Omega^n B_n \Xtilde_c(T_1)$, and we indeed recover the statement of
  Theorem~\ref{MAIN}.
\end{proof}

\begin{proof}[Proof of Corollary~$\ref{APX}$]
  By \cite{schwanzl-vogt} the free $F_n$-algebra generated by a
  CW-complex $X$ is cofibrant in $\Alg^{\Tn}$. The space $B_n F_n X$
  is fibrant, as any object of $R_{S^n}\Spaces$. Then the isomorphism
  $B_n F_n X \xrightarrow{\id} B_n F_n X$ implies by
  Theorem~\ref{MAIN2} that the adjoint $F_n X \ra \Omega_n B_n F_n X$
  is a weak equivalence. On the other hand, $B_n F_n X = \Sigma^n X$
  by \ref{BF}.  Thus, we get a weak equivalence $F_n X
  \xrightarrow{\sim} \Omega^n \Sigma^n X$. It defines an equivalence
  of monads, because of the naturality of the construction.
\end{proof}

\begin{proof}[Proof of Theorem~$\ref{MAIN2}$]
  It is enough to show that for every cofibrant $\Tn$-algebra
  $\Xtilde$, the unit $\eta_{\Xtilde} \colon \Xtilde\ra
  \Omega^nB_n\Xtilde$ of the adjunction $(B_n, \Omega^n)$ is a weak
  equivalence in $\Alg^{\Tn}$. Indeed, for $\Xtilde\in \Alg^{\Tn}$,
  $Y\in\Spaces$, and $f\colon \Xtilde\ra \Omega^nY$, we have a
  commutative diagram
  $$\xymatrix{ \Xtilde\ar[dr]_-f\ar[r]^-{\eta_{\Xtilde}}&
    \Omega^nB_n\Xtilde
    \ar[d]^{{\Omega^nf^\flat}}\\
    &\Omega^nY, }$$
  where $f^\flat$ is the adjoint to $f$.  Assume
  that $\Xtilde$ is cofibrant. By assumption $\eta_{\Xtilde}$ is a
  weak equivalence in $\Alg^{\Tn}$. If $f$ is also a weak equivalence,
  then so is $\Omega^nf^\flat$.  In particular the map
  $$\Omega^nf^\flat(T_1)\colon \Omega^n (B_n \Xtilde) = (\Omega^n B_n
  \Xtilde) (T_1) \ra (\Omega^n Y)(T_1) = \Omega^n Y$$
  is a weak
  equivalence of spaces, or, in other words, $f^\flat$ is an
  $S^n$-local weak equivalence.
  
  Conversely, if $f^\flat$ is an $S^n$-local equivalence, then
  $\Omega^nf^\flat$ is an objectwise weak equivalence, and so is $f$.
  
  The proof  of the fact that for a cofibrant $\Xtilde \in\Alg^{\Tn}$,
  the map $\eta_{\Xtilde}$ is a weak equivalence follows from a
  bootstrap argument below.
\smallskip

\noindent
1) Let $\Xtilde = F_n(Z)$, where $Z$ is an arbitrary pointed discrete
space.  Since $F_n$ is a left adjoint functor, it commutes with
colimits. Therefore, since $Z$ is the colimit of the poset of finite
subsets $Y$ of $Z$ containing the basepoint, we get:
$$F_n(Z) = \colim_{Y\subseteq Z} F_n(Y) = \colim_{Y\subseteq Z}
\Map(S^n, Y\wedge S^n).$$
The second
equality follows from \ref{KEY}.  Furthermore, since $S^n$ is a
compact space, we have $ \colim_{Y\subseteq Z} \Map(S^n, Y\wedge
S^n) \linebreak[2] = \linebreak[1] \Map(S^n, \linebreak[0] Z\wedge
S^n)$.  Therefore, the map $\eta_{\Xtilde}$ is an isomorphism of
$\Tn$-algebras by \ref{BF}.  \smallskip

\noindent
2) Let $Z_\bullet$ be a pointed simplicial set, and let $\Xtilde= F_n'
(Z_\bullet)$, where $F'_n = F'_{\Tn}$.  We have by \ref{FREE ALG}
$$F'_n (Z_\bullet ) = F_n(|Z_\bullet |) \cong | F_n Z_\bullet |,$$
where $F_n Z_\bullet$ denotes the simplicial $\Tn$-algebra obtained by
applying $F_n$ in each  simplicial dimension of $Z_\bullet$.  By Step 1
for every $k\geq 0 $, we have an isomorphism $\eta_k \colon F_n(Z_k)
\ra \Omega^nB_nF_n(Z_k)$, assembling into a simplicial map by
naturality. Thus, the map
$$|\eta_\bullet |\colon \Xtilde \ra |\Omega^nB_nF_n(Z_\bullet)|$$
is
also an isomorphism.  Next,  notice that by \ref{BF}, we have
$B_nF_n(Z_k) = {Z_k}\wedge S^n$, so it is an $(n-1)$-connected
space.  Therefore (see \cite[Theorem 12.3]{May}), we have a natural
weak equivalence $|\Omega^nB_nF_n(Z_\bullet)|\simeq \Omega^n
|B_nF_n(Z_\bullet)|$. (A technical condition of properness of $B_n F_n
(Z_\bullet)$, needed for applying May's theorem, is satisfied here, as
$Z_\bullet$ is discrete and $B_n$ and $F_n$ are admissible functors,
see \cite[Definitions 11.2 and A.7]{May}.)  Combining this with the
isomorphism $\ |B_nF_n(Z_\bullet)| \cong B_n|F_n(Z_\bullet)|$ we get a
weak equivalence
$$|\Omega^nB_nF_n(Z_\bullet)|\simeq \Omega^nB_n \abs{F_n (Z_\bullet)}
\cong \Omega^nB_n\Xtilde$$
It follows that $\eta_{\Xtilde}$ is a weak
equivalence.
\smallskip
  
\noindent
3) Let $\Xtilde$ be any  $\Tn$-algebra and
${F'_{n}U'_{n}}_\bullet \Xtilde$ its simplicial resolution as in
Section~\ref{FUA}, where $U'_n = U'_{\Tn}$.  Note that, in every
simplicial dimension $k$, the algebra $(F'_{n}U'_{n})_k \Xtilde$ is of
the form considered in Step 2. It follows that for $k\geq 0$, we have
a weak equivalence
\begin{equation}
\label{spaces}
\eta_k \colon (F'_{n}U'_{n})_k \Xtilde \xrightarrow{\sim} \Omega^nB_n
  (F'_{n}U'_{n})_k \Xtilde.
\end{equation}  

To see that the map
$$|\eta_\bullet | \colon |{F'_{n}U'_{n}}_\bullet \Xtilde| \ra
|\Omega^nB_n {F'_{n}U'_{n}}_\bullet \Xtilde|$$
is also a weak
equivalence, we can use a result of May \cite[Theorem 11.13]{May}.
The assumption of strict properness \cite[Definition 11.2]{May} of the
simplicial spaces ${F'_n U'_n}_\bullet \Xtilde$ and $\Omega^n B_n
{F'_n U'_n}_\bullet \Xtilde$, needed for May's theorem, is not hard to
verify, since all the functors $F_n$, $U_n$, $\abs{\Sing (\cdot)}$,
$B_n$, and $\Omega_n$ are admissible in the sense of \cite[Definition
A.7]{May}. May also assumes that the realizations of the simplicial
spaces are connected H-spaces, which will not be satisfied in our
case, in general. His result however readily generalizes to the case
of simplicial spaces whose realizations are H-spaces with $\pi_0$'s
having a group structure, as it is the case for the simplicial spaces
at hand for $n \ge 1$. The H-space structure is not there for $n=0$,
but in this case, the statement of the theorem is trivial, anyway.
  
Using arguments similar to those employed in Step 2, we get from here
that
$$\eta \colon |{F'_{n}U'_{n}}_\bullet \Xtilde| \ra \Omega^nB_n
|{F'_{n}U'_{n}}_\bullet \Xtilde|$$
is a weak equivalence. 

\noindent 
4) Let $\Xtilde$ be any cofibrant algebra. We have
a commutative diagram:
$$\xymatrix{ \abs{{F'_n U'_n}_\bullet \Xtilde} \ar[d]^{\sim}_h
  \ar[r]^{\eta}_{\sim} &
  \Omega^n B_n \abs{{F'_n U'_n}_\bullet \Xtilde} \ar[d]^{\Omega^n B_n h}\\
  \Xtilde \ar[r]^{\eta_{\Xtilde}} & \Omega^n B_n \Xtilde , }$$
where
$h$ is the weak equivalence of Lemma~\ref{FUA WE}. The functor $B_n$
is a left Quillen functor and as such it preserves weak equivalences between
cofibrant $\Tn$-algebras, while  $\Omega^n$ preserves all weak
equivalences. Therefore  $\Omega^n B_n h$ is a weak equivalence, and, as a
consequence, so is $\eta_{\Xtilde}$.
\end{proof}

\begin{thm}
\label{unique}
Suppose $\TT$ is an algebraic theory such that it
\begin{enumerate}
\item acts on $n$-fold loops spaces $\Omega^n X$ by natural operations
  $(\Omega^n X)^k \to (\Omega^n X)^l$, i.e., admits a morphism $\phi:
  \TT \to \Tn$, and
\item via this action deloops $n$-fold loop spaces in the sense of
  Theorem~$\ref{MAIN2}$, i.e., the loop functor $R_A \Spaces
  \xrightarrow{\Omega^n} \Alg^{\Tn} \xrightarrow{\phi^*} \Alg^{\TT}$
  establishes a Quillen equivalence.
\end{enumerate}
Then $\phi: \TT \to \Tn$ is a weak equivalence of topological theories.
\end{thm}
This theorem is, in fact, an obvious corollary of a uniqueness theorem
\cite[Theorem 1.6]{badzioch:gamma} (theories considered in 
\cite{badzioch:gamma} are enriched over simplicial sets, but the proof of this
result holds for topological theories with little changes).


\providecommand{\bysame}{\leavevmode\hbox to3em{\hrulefill}\thinspace}
\providecommand{\MR}{\relax\ifhmode\unskip\space\fi MR }
\providecommand{\MRhref}[2]{%
  \href{http://www.ams.org/mathscinet-getitem?mr=#1}{#2}
}
\providecommand{\href}[2]{#2}

\end{document}

%% file: loopmacro.tex
\newcommand{\CC}{\mathcal{C}}
\newcommand{\DD}{\mathcal{D}}
\newcommand{\OO}{\mathcal{O}}
\newcommand{\TT}{\mathcal{T}}

\newcommand{\Alg}{\mathcal{A}lg}
\newcommand{\Tn}{\TT^{S^n}}

\newcommand{\abs}[1]{\lvert#1\rvert}
\newcommand{\End}[1]{\operatorname{End}(#1)}

\newcommand{\Map}{\operatorname{Map}_*}

\newcommand{\Mor}{\operatorname{Mor}}
\newcommand{\op}{{\operatorname{op}}}

\newcommand{\Xtilde}{\widetilde{X}}

\newcommand{\id}{\operatorname{id}}
\newcommand{\ev}{\operatorname{ev}}

\newcommand{\colim}{\operatorname{colim}}

\newcommand{\Sing}{\operatorname{Sing}_\bullet}
\newcommand{\Spaces}{\operatorname{Spaces}_*}
\newcommand{\SSets}{\operatorname{SSets}_*}

\newcommand{\ra}{\rightarrow}

\newtheorem{thm}{Theorem}[section]
\newtheorem{lm}[thm]{Lemma}
\newtheorem{prop}[thm]{Proposition}
\newtheorem{crl}[thm]{Corollary}

\theoremstyle{definition}

\newtheorem{df}[thm]{Definition}
\newtheorem{ex}[thm]{Example}
\newtheorem{notation}[thm]{Notation}

\newtheorem{nn}[thm]{}

%% file: bcv.bbl
\begin{thebibliography}{May72}

\bibitem[Bad02]{ATHT}
B.~Badzioch, \emph{Algebraic theories in homotopy theory}, Ann. of Math. (2)
  \textbf{155} (2002), no.~3, 895--913. \MR{2003g:55035}

\bibitem[Bad03]{badzioch:gamma}
\bysame, \emph{{From $\Gamma$-spaces to algebraic theories}}, Preprint,
  University of Minnesota, June 2003, \texttt{math.AT/0306010}.

\bibitem[Bou92]{Bousfield}
A.~K. Bousfield, \emph{The simplicial homotopy theory of iterated loop spaces},
  Manuscript, 1992.

\bibitem[BV73]{boardman-vogt}
J.~M. Boardman and R.~M. Vogt, \emph{Homotopy invariant algebraic structures on
  topological spaces}, Springer-Verlag, Berlin, 1973, Lecture Notes in
  Mathematics, Vol. 347. \MR{54 \#8623a}

\bibitem[Chu04]{chung}
K.~Chung, \emph{{Ph.D. thesis}}, In progress, University of Minnesota, 2004.

\bibitem[Coh76]{cohen:cubes}
F.~R. Cohen, \emph{The homology of {$\mathcal{C}_{n+1}$}-spaces, {$n\ge0$}},
  The homology of iterated loop spaces, Lecture Notes in Math., vol. 533,
  Springer-Verlag, 1976, pp.~207--351.

\bibitem[Coh88]{cohen:sundry}
\bysame, \emph{Artin's braid groups, classical homotopy theory and sundry other
  curiosities}, Contemp. Math. \textbf{78} (1988), 167--206.

\bibitem[CS99]{chas-sullivan}
M.~Chas and D.~Sullivan, \emph{String topology}, Preprint, CUNY, November 1999,
  \texttt{math.GT/9911159}.

\bibitem[Hir03]{hirschhorn}
P.~S. Hirschhorn, \emph{Model categories and their localizations}, Mathematical
  Surveys and Monographs, vol.~99, American Mathematical Society, Providence,
  RI, 2003. \MR{2003j:18018}

\bibitem[Hov99]{hovey}
M.~Hovey, \emph{Model categories}, Mathematical Surveys and Monographs,
  vol.~63, American Mathematical Society, Providence, RI, 1999. \MR{99h:55031}

\bibitem[May72]{May}
J.~P. May, \emph{The geometry of iterated loop spaces}, Springer-Verlag,
  Berlin, 1972, Lectures Notes in Mathematics, Vol. 271. \MR{54 \#8623b}

\bibitem[ML98]{ML}
S.~Mac~Lane, \emph{Categories for the working mathematician}, second ed.,
  Graduate Texts in Mathematics, vol.~5, Springer-Verlag, New York, 1998.
  \MR{2001j:18001}

\bibitem[MT78]{may-thomason}
J.~P. May and R.~Thomason, \emph{The uniqueness of infinite loop space
  machines}, Topology \textbf{17} (1978), no.~3, 205--224. \MR{80g:55015}

\bibitem[Sal03]{salvatore}
P.~Salvatore, \emph{The universal operad of iterated loop spaces}, Draft, 2003.

\bibitem[Seg74]{segal}
G.~Segal, \emph{Operations in stable homotopy theory}, New developments in
  topology (Proc. Sympos. Algebraic Topology, Oxford, 1972), Cambridge Univ.
  Press, London, 1974, pp.~105--110. London Math Soc. Lecture Note Ser., No.
  11. \MR{49 \#3917}

\bibitem[Sta63]{stasheff}
J.~Stasheff, \emph{Homotopy associativity of {$H$}-spaces. {I}, {II}}, Trans.
  Amer. Math. Soc. \textbf{108} (1963), 275-292; ibid. \textbf{108} (1963),
  293--312. \MR{28 \#1623}

\bibitem[SV91]{schwanzl-vogt}
R.~Schw{\"a}nzl and R.~M. Vogt, \emph{The categories of {$A\sb \infty$}- and
  {$E\sb \infty$}-monoids and ring spaces as closed simplicial and topological
  model categories}, Arch. Math. (Basel) \textbf{56} (1991), no.~4, 405--411.
  \MR{92b:18006}

\bibitem[SV04]{sullivan-v}
D.~Sullivan and A.~A. Voronov, \emph{Brane topology}, Draft, 2004.

\bibitem[Tho79]{thomason}
R.~W. Thomason, \emph{Uniqueness of delooping machines}, Duke Math. J.
  \textbf{46} (1979), no.~2, 217--252. \MR{80e:55013}

\end{thebibliography}
